\documentclass[11pt,twoside,dvips,letterpaper]{paper}

\usepackage{graphicx}
\usepackage{amsmath}
\usepackage{amssymb}
\usepackage{dcolumn}
\usepackage{multirow}
\usepackage{bm}
\usepackage[usenames,dvipsnames]{color}
\usepackage{alltt}
\usepackage[colorlinks,dvips]{hyperref}
\usepackage[hyphenbreaks]{breakurl}
\usepackage[margin=1in]{geometry}
\usepackage[compact]{titlesec}
\usepackage{mdwlist}

\newcommand{\bv}[1]{{\mbox{\boldmath$\mathbf{#1}$}}}
\newcommand{\rmean}[1]{{\left\langle{#1}\right\rangle}}
\newcommand{\irmean}[1]{{\langle{#1}\rangle}}
\newcommand{\irv}[1]{\langle{#1^2}\rangle}
\newcommand{\rs}[1]{\rmean{#1^3}}
\newcommand{\irs}[1]{\langle{#1^3}\rangle}
\newcommand{\rk}[1]{\rmean{#1^4}}
\newcommand{\irk}[1]{\langle{#1^4}\rangle}

\newcommand{\Eqr}[1]{(\ref{#1})}
\newcommand{\Eqre}[1]{Eq.~(\ref{#1})}

\newcommand{\Eqres}[1]{Eqs.~(\ref{#1})}

\newcommand{\Fige}[1]{Figure~\ref{#1}}

\newcommand{\laur}{LA-UR-13-28548}

\newcommand{\ti}{Diffusion processes satisfying a conservation law constraint}
\newcommand{\aufirst}{J.\ Bakosi}
\newcommand{\ausecond}{J.R.\ Ristorcelli}
\newcommand{\au}{\aufirst and \ausecond}
\newcommand{\kw}{Conservation;
                 Diffusion process;
                 Fokker-Planck equation;
                 Statistical moment equations}

\hypersetup{citecolor=blue,
            linkcolor=blue,
            urlcolor=blue,
            pdftitle=\ti,
            pdfauthor=\au,
            pdfkeywords=\kw}
\usepackage{hypernat}

\begin{document}

\title{\ti\\[0.2cm]\small\texttt{\laur}\\[0.2cm]\texttt{Accepted for publication in
International Journal of Stochastic Analysis, January 5, 2014}\\[-0.3cm]}

\author{\normalsize\aufirst, \normalsize\ausecond\\[-0.1cm]
        \texttt{\normalsize\{jbakosi,jrrj\}@lanl.gov} \\[-0.1cm]
        \normalsize Los Alamos National Laboratory, Los Alamos, NM 87545,
        USA\\[-0.3cm]}

\maketitle

\begin{abstract}
We investigate coupled stochastic differential equations governing $N$
non-negative continuous random variables that satisfy a conservation principle.
In various fields a conservation law requires that a set of fluctuating
variables be non-negative and (if appropriately normalized) sum to one. As a
result, any stochastic differential equation model to be realizable must not
produce events outside of the allowed sample space. We develop a set of
constraints on the drift and diffusion terms of such stochastic models to ensure
that both the non-negativity and the unit-sum conservation law constraint are
satisfied as the variables evolve in time. We investigate the consequences of
the developed constraints on the Fokker-Planck equation, the associated system
of stochastic differential equations, and the evolution equations of the first
four moments of the probability density function. We show that random variables,
satisfying a conservation law constraint, represented by stochastic diffusion
processes, must have diffusion terms that are coupled and nonlinear. The set of
constraints developed enables the development of statistical representations of
fluctuating variables satisfying a conservation law. We exemplify the results
with the bivariate beta process and the multivariate Wright-Fisher, Dirichlet,
and Lochner's generalized Dirichlet processes.
\keywords{\kw}
\end{abstract}

\section{Introduction and problem statement}
We investigate the consequences of the unit-sum requirement on $N\!>\!1$
non-negative continuous random variables governed by a diffusion process. Such
mathematical description is useful to represent fluctuating variables,
$Y_1,\dots,Y_N$, subject to the constraint $\sum Y_\alpha\!=\!1$. We are
interested in stochastic diffusion models and statistical moment equations
describing the temporal evolutions $Y_\alpha=Y_\alpha(t)$ and their statistics.
In particular, we study the consequences of the bounded sample space, required
by the non-negativity of $Y_\alpha$ and the unit-sum conservation principle,
$\sum Y_\alpha\!=\!1$. A simple physical example is the mixture of different
chemical species, represented by mass fractions $0 \le Y_\alpha \le 1$
undergoing reaction in a fluid whose overall mass is conserved. Such
mathematical problems also appear in evolutionary theory \cite{Pearson_1896},
Bayesian statistics \cite{Paulino_95}, geology \cite{Chayes_62, Chayes_66,
Martin_65}, forensics \cite{Lange_95}, econometrics \cite{Gourieroux_06},
turbulent mixing and combustion \cite{Girimaji_91}, and population biology
\cite{Steinrucken_2013}. Mathematical properties of such random fractions are
given in \cite{Mauldon_51, Mauldon_59, Good_65, Pyke_65}.

Mathematically, we are interested in the following question: What functions are
allowed to represent the drift, $A_\alpha$, and diffusion, $b_{\alpha\beta}$,
terms of the system, governing the vector $\bv{Y}=(Y_1,\dots,Y_N)$,
\begin{equation}
\mathrm{d}Y_\alpha(t) = A_\alpha(\bv{Y},t)\mathrm{d}t + \sum_{\beta=1}^N
b_{\alpha\beta}(\bv{Y},t)\mathrm{d}W_\beta(t), \qquad \alpha = 1,\dots,N,
\label{eq:Ito}
\end{equation}
if
\begin{equation}
Y_\alpha \ge 0, \quad \alpha=1,\dots,N \qquad \textrm{and} \qquad
\sum_{\alpha=1}^NY_\alpha=1, \label{eq:realizable}
\end{equation}
must hold for all $t$. In \Eqre{eq:Ito} $\mathrm{d}W_\alpha(t)$ is a
vector-valued Wiener process with mean $\irmean{\mathrm{d}W_\alpha} = 0$ and
covariance $\irmean{\mathrm{d}W_\alpha \mathrm{d}W_\beta} = \delta_{\alpha\beta}
\mathrm{d}t$, see \cite{Gardiner_09}, and $\delta_{\alpha\beta}$ is Kronecker's
delta. If the components of $\bv{Y}$ satisfy the constraints in
\Eqre{eq:realizable}, we call the event $\bv{Y}$ realizable. A consequence of
the constraints in \Eqre{eq:realizable} imposed on the stochastic system
\Eqre{eq:Ito} is that for all $t$ the following holds
\begin{equation}
\sum_{\alpha=1}^N \mathrm{d}Y_\alpha(t) = 0 = \sum_{\alpha=1}^N
A_\alpha(\bv{Y},t)\mathrm{d}t + \sum_{\alpha=1}^N\sum_{\beta=1}^N
b_{\alpha\beta}(\bv{Y},t)\mathrm{d}W_\beta(t).
\label{eq:sumIto}
\end{equation}
In other words, we are interested in expressions for $A_\alpha$ and
$b_{\alpha\beta}$, what constraints they must satisfy in addition to
\Eqre{eq:sumIto}, and how to implement them so that \Eqre{eq:Ito} produces
realizable events, i.e., $\bv{Y}$ satisfies \Eqre{eq:realizable} for all $t$.

We study diffusion processes as: (1) they are mathematically simple vehicles for
representing temporal evolutions of fluctuating fractions (of a
unit) and their statistics, (2) they lend themselves to simple
Monte-Carlo numerical methods \cite{Kloeden_99}, and (3) they serve as a
starting point for representations of statistical moment equations if individual
samples and joint probabilities are not required.  The Markovian assumption
\cite{Gardiner_09} is made at the outset and jump contributions are ignored. We
derive constraints for the drift and diffusion terms that assure that the
modeled processes are realizable (i.e., produce non-negative variables that
satisfy the unit-sum constraint) for any realization at all times. We address
the problem of the functional forms of the drift and diffusion terms from three
perspectives: (1) the Fokker-Planck equation for the probability density
function, (2) the stochastic differential equations for the individual
realizations, and (3) the evolution equations for the jointly coupled
statistics.

The plan of the paper is as follows. \S\ref{sec:samplespace} introduces the
geometry of the multi-dimensional sample space within which realizations of
fractions of a unit are allowed and discusses constraints that ensure realizable
statistical moments. \S\ref{sec:PDF} develops the implications of realizability
on diffusion processes governing fractions. \S\ref{sec:momentevolution} follows
by developing realizability constraints on the time-evolutions of statistics.
\S\ref{sec:models} surveys some existing realizable diffusion processes. A
summary is given in \S\ref{sec:summary}.

\section{Realizability due to conservation}
\label{sec:samplespace}
The notion of realizability due to a conservation law constraint was introduced
and defined by \Eqre{eq:realizable}. We now discuss the consequences of
realizability pertaining to the individual samples of the state space,
\S\ref{sec:realsamples}, and of their statistics, \S\ref{sec:realmoments}.

\subsection{The universal geometry of allowed realizations}
\label{sec:realsamples}
The geometrical definition of the sample space is given in which the vector
$\bv{Y}=(Y_1,\dots,Y_N)$ is allowed if \Eqre{eq:realizable} is to be satisfied.
This is used to derive constraints for stochastic diffusions and their moment
equations in the subsequent sections.

A realization of the vector, $\bv{Y}$, with coordinates $Y_\alpha \ge 0$,
$\alpha=1,\dots,N$, specifies a point in the multi-dimensional sample space. The
union of all such points that satisfy
\begin{equation}
\sum_{\alpha=1}^NY_\alpha = 1,
\label{eq:unitsum}
\end{equation}
is the space of allowed realizations, \Eqre{eq:realizable}. For example, in
representing mass fraction constituents of a substance, \Eqre{eq:unitsum}
restricts the possible components of $\bv{Y}$ to those that are realizable;
those vectors that point outside of the allowed space are not conserved; if
\Eqre{eq:unitsum} is violated, spurious mass is created or destroyed.

\begin{figure}
\centering
{\resizebox{0.4\columnwidth}{!}{\input{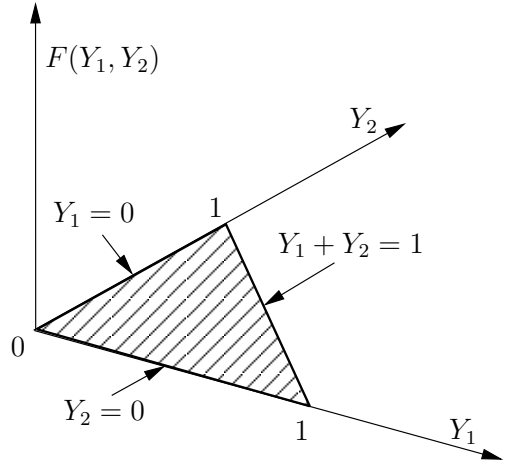}}}
\caption{The geometry of allowed realizations for $N=3$ variables, $\bv{Y} =
(Y_1, Y_2, Y_3 = 1 - Y_1 - Y_2)$, that satisfy non-negativity and the unit-sum
constraint, \Eqre{eq:realizable}. The boundary of the allowed region on the
plane spanned by $Y_1$ and $Y_2$ is the closed loop of straight lines, $(Y_1 =
0, \enskip Y_2 = 0, \enskip Y_1 + Y_2 = 1)$, defined by \Eqre{eq:boundary}. If
the vector $\bv{Y}$ points inside the triangle, it is a realizable event.}
\label{fig:simplex}
\end{figure}

Mathematically, the geometry of allowed realizations is a simplex, the
generalization of a triangle to multiple dimensions. For $N$ variables the
$(N\!-\!1)$-simplex is a bounded convex polytope, $\mathcal{P}$, on the
$(N\!-\!1)$-dimensional hyperplane; $\mathcal{P}$ is the convex hull of its $N$
vertices. $\mathcal{P}$'s boundary, $\partial\mathcal{P}$, is defined as the
closed surface of non-overlapping hyperplanes of $N-2$ dimensions,
\begin{equation}
\partial\mathcal{P} \equiv \left( Y_\alpha=0 \enskip : \enskip
\alpha=1,\dots,N-1; \quad \sum_{\alpha=1}^{N-1} Y_\alpha = 1\right),
\label{eq:boundary}
\end{equation}
plotted in \Fige{fig:simplex} for $N=3$.

The domain (or support) of the joint probability, $F(\bv{Y})$, with
$Y_\alpha$, $\alpha=1,\dots,N$, is the $(N\!-\!1)$-simplex. Of all $Y_\alpha$
only $N\!-\!1$ are independent due to \Eqre{eq:unitsum} and without loss of
generality we take
\begin{equation}
Y_N = 1-\sum_{\alpha=1}^{N-1}Y_\alpha.
\label{eq:YN}
\end{equation}
The same geometry of allowed realizations is discussed by Pope in the
$N$-dimensional state space in \cite{Pope_04mix} in the context of ideal gas
mixing in turbulent combustion.

We confine our attention here to $N\!-\!1$ dimensions, as one of the variables
is determined by the unit-sum requirement, \Eqre{eq:YN}. As a consequence, the
$(N\!-\!1)$-dimensional geometry of realizable events is remarkably simple and
universal: it is the bounded convex polytope whose boundary is defined by
\Eqre{eq:boundary}. Consequently, the realizability constraint,
\Eqre{eq:realizable}, uniquely and universally determines the realizable region
of the state space: it is the same in all points in space and time for all
materials undergoing any physical process that conserves mass,
\Eqre{eq:unitsum}. The ensemble is realizable if and only if all samples reside
inside the polytope given by \Eqre{eq:boundary}. For $N\!=\!3$ this means that
the support of $F$ is the triangle depicted in \Fige{fig:simplex}.

\subsection{Realizable statistical moments}
\label{sec:realmoments}
If the fractions are non-negative and sum to one, \Eqre{eq:realizable}, they are
also bounded,
\begin{equation}
0 \le Y_\alpha \le 1, \enskip \alpha=1,\dots,N, \qquad \textrm{and} \qquad
\sum_{\alpha=1}^NY_\alpha=1, \label{eq:bounded}
\end{equation}
whose consequences on some of their statistical moments are now discussed.

Taking mathematical expectations of \Eqre{eq:bounded}, see e.g.\
\cite{vanKampen_04}, yields
\begin{equation}
0 \le \irmean{Y_\alpha} \le 1, \enskip \alpha=1,\dots,N \qquad \textrm{and}
\qquad \sum_{\alpha=1}^N\irmean{Y_\alpha}=1. \label{eq:boundedmeans}
\end{equation}
Similar to the instantaneous fractions, the first statistical moments are also
non-negative, bounded, and sum to unity.

Since both the instantaneous variables and their means are bounded,
fluctuations about the means are also bounded:
\begin{equation}
-1 \le y_\alpha = Y_\alpha - \irmean{Y_\alpha} \le 1.
\label{eq:boundedfluctuations}
\end{equation}
As a consequence, the variances and the covariances are also bounded:
\begin{align}
0 \le \irv{y_\alpha} & = \rmean{(Y_\alpha - \irmean{Y_\alpha})^2} \le 1,
\label{eq:boundedvariances} \\
-1 \le \irmean{y_\alpha y_\beta} & = \rmean{(Y_\alpha -
\irmean{Y_\alpha})(Y_\beta - \irmean{Y_\beta})} \le 1, \qquad \alpha\ne\beta.
\label{eq:boundedcovariances}
\end{align}
Multiplying \Eqre{eq:unitsum} by $y_\beta$, $\beta\!=\!1,\dots,N$ and taking
the expectation yield
\begin{align}
\irv{y_1} + \irmean{y_2y_1} + \dots + \irmean{y_Ny_1} & = 0 \nonumber\\
\irmean{y_1y_2} + \irv{y_2} + \dots + \irmean{y_Ny_2} & = 0 \nonumber\\
& \vdots \label{eq:zerosums}\\
\irmean{y_1y_N} + \irmean{y_2y_N} + \dots + \irv{y_N} & = 0, \nonumber
\end{align}
i.e., the row-sums and, due to symmetry, the column-sums of the covariance
matrix are zero. Expressing $\irmean{y_Ny_1}$, $\irmean{y_Ny_2}$, etc., from
the first $N\!-\!1$ equations of \Eqres{eq:zerosums} and substituting them into
the $N^\mathrm{th}$ one, yield the weaker constraint,
\begin{equation}
\sum_{\alpha=1}^{N-1} \sum_{\beta=1}^{N-1} \irmean{y_\alpha y_\beta} -
\irv{y_N} = 0. \label{eq:weakzerosums}
\end{equation}
Due to bounded fluctuations, \Eqre{eq:boundedfluctuations}, the third central
moments are also bounded,
\begin{equation}
-1 \le \irs{y_\alpha} = \rmean{(Y_\alpha - \irmean{Y_\alpha})^3} \le 1,
\label{eq:boundedthirdmoments}
\end{equation}
and in general, for $n \ge 2$ we have,
\begin{align}
0 \le \irmean{y_\alpha^n} \le 1, & \qquad \textrm{for even $n$},
\label{eq:boundedneven} \\
-1 \le \irmean{y_\alpha^n} \le 1, & \qquad \textrm{for odd $n$}.
\label{eq:boundednodd}
\end{align}

Ensuring non-negativity and unit-sum puts constraints on possible time
evolutions of $\bv{Y}\!=\!\bv{Y}(t)$, represented by diffusion processes and
that of their statistics. Some of these constraints are developed in the
following sections.

\section{Diffusion processes for random fractions}
\label{sec:PDF}
Implications of the geometry of the realizable state space, discussed in
\S\ref{sec:samplespace}, on diffusion processes are developed. First, the
relevant mathematical properties of Fokker-Planck equations are reviewed in
\S\ref{sec:math}, followed by the constraints on their functional forms,
\S\ref{sec:sdes}.

\subsection{Review of some boundary conditions of Fokker-Planck equations}
\label{sec:math}
The discussion is restricted to Markov processes which by definition obey a
Chapman-Kolmogorov equation \cite{Gardiner_09}. Assuming that $Y_\alpha$ are
continuous in space and time, jump processes are excluded. The temporal
evolution of random fractions, $\bv{Y}(t)$, constrained by \Eqre{eq:realizable}
can then be represented most generally by diffusion processes whose transitional
probability, $F(\bv{Y},t)$, is governed by the Fokker-Planck equation,
\begin{equation}
\frac{\partial}{\partial t}F(\bv{Y},t) = - \sum_{\alpha=1}^{N-1}
\frac{\partial}{\partial Y_\alpha}\big[ A_\alpha(\bv{Y},t) F(\bv{Y},t)
\big] + \frac{1}{2} \sum_{\alpha=1}^{N-1} \sum_{\beta=1}^{N-1}
\frac{\partial^2}{\partial Y_\alpha \partial
Y_\beta}\big[B_{\alpha\beta}(\bv{Y},t)F(\bv{Y},t)\big],\label{eq:FP}
\end{equation}
where $A_\alpha$ and $B_{\alpha\beta}$ denote drift and diffusion in state
space, respectively, and $B_{\alpha\beta}$ is symmetric non-negative
semi-definite \cite{vanKampen_04}. \Eqre{eq:FP} is a partial differential
equation that governs the joint probability, $F(\bv{Y},t)$, of the fractions,
$Y_\alpha$, $\alpha = 1,\dots,N-1$. $Y_N$ is excluded from \Eqre{eq:FP} and is
determined by \Eqre{eq:YN}. Augmented by initial and boundary conditions,
\Eqre{eq:FP} describes the transport of probability in sample space
$\mathcal{R}$ whose boundary is $\partial\mathcal{R}$ with normal vector
$n_\alpha$ , see \cite{Gardiner_09}.

\Eqre{eq:FP} can be written in conservation form as
\begin{equation}
\frac{\partial}{\partial t} F(\bv{Y},t) + \sum_{\alpha=1}^{N-1}
\frac{\partial}{\partial Y_\alpha} I_\alpha(\bv{Y},t) = 0, \label{eq:FPI}
\end{equation}
in terms of the probability flux, see \cite{Gardiner_09}, Sec.\ 5.1,
\begin{equation}
\begin{split}
I_\alpha(\bv{Y},t) = A_\alpha(\bv{Y},t) F (\bv{Y},t) - \frac{1}{2}
\sum_{\beta=1}^{N-1} \frac{\partial}{\partial Y_\beta}
\big[B_{\alpha\beta}(\bv{Y},t) F(\bv{Y},t)\big], \label{eq:probflux}\\
\qquad\qquad \alpha=1,\dots,N-1.
\end{split}
\end{equation}
Using Eqs.\ (\ref{eq:FPI}--\ref{eq:probflux}) the following boundary conditions
are considered, see \cite{Gardiner_09}, Sec.\ 6.2.

\begin{enumerate}
\item \emph{Reflecting barrier.} If $\sum n_\alpha I_\alpha(\bv{Y},t) = 0$
everywhere on the boundary, $\partial\mathcal{R}$ is a reflecting barrier: a
particle inside $\mathcal{R}$ cannot cross the boundary and must be reflected
there.
\item \emph{Absorbing barrier.} If $F(\bv{Y},t) = 0$ everywhere on
the boundary, $\partial\mathcal{R}$ is an absorbing barrier: if a particle
reaches the boundary, it is removed from the system.
\item \emph{Other types of boundary conditions.} Some part of the boundary may
be reflecting while some other may be absorbing: a combination is certainly
possible. We only consider reflecting and absorbing barriers -- other types of
boundaries are discussed in \cite{Feller_52}.
\end{enumerate}

To support the forthcoming discussion, some well-established mathematical
properties of multi-variable Fokker-Planck equations have been reviewed.

\subsection{Realizable diffusion processes}
\label{sec:sdes}
The implications of the realizability constraint, \Eqre{eq:realizable}, on the
functional forms of the drift and diffusion terms of the Fokker-Planck equation
\Eqr{eq:FP} are now investigated.

As discussed in \S\ref{sec:samplespace}, the region of the sample space allowed
by the realizbility requirement is the polytope $\mathcal{P}$ defined by its
boundary, $\partial\mathcal{P}$, \Eqre{eq:boundary}, in which all samples of
$\bv{Y}\!=\!\bv{Y}(t)$ must reside at all times. Consequently, the sample space,
$\mathcal{R}$, of the Fokker-Planck equation \Eqr{eq:FP} must coincide with
$\mathcal{P}$, which constrains the possible functional forms of
$A_\alpha(\bv{Y},t)$ and $B_{\alpha\beta}(\bv{Y},t)$. In the following, these
constraints are developed for binary (single-variable) processes first, followed
by ternary processes, and then generalized to multiple variables.

\subsubsection{Realizable binary processes: $\bv{N=2}$}
The It\^o diffusion process \cite{Gardiner_09}, governing the variable $Y$,
\begin{equation}
\mathrm{d}Y(t) = A(Y,t)\mathrm{d}t + \sqrt{B(Y,t)}\mathrm{d}W(t),
\label{eq:sIto}
\end{equation}
with $B(Y,t)\ge0$ is equivalent to and derived from \Eqre{eq:FP} with $N=2$,
see e.g., \cite{Gardiner_09},
\begin{equation}
\frac{\partial}{\partial t}F(Y,t) = -\frac{\partial}{\partial
Y}\big[A(Y,t)F(Y,t)\big] + \frac{1}{2}\frac{\partial^2}{\partial
Y^2}\big[B(Y,t)F(Y,t)\big]. \label{eq:sFP}
\end{equation}
For $N=2$ the allowed space of realizations is a line with endpoints given by
\Eqre{eq:boundary},
\begin{equation}
\big( Y = 0; \quad Y = 1 \big). \label{eq:endpoints}
\end{equation}
This can be ensured if the drift and diffusion terms in Eqs.\
(\ref{eq:sIto}--\ref{eq:sFP}) satisfy
\begin{align}
A(Y\!=\!0,t) \ge 0 \qquad \mathrm{and} \qquad B(Y\!=\!0,t) = 0,\label{eq:scs1}\\
A(Y\!=\!1,t) \le 0 \qquad \mathrm{and} \qquad B(Y\!=\!1,t) = 0.\label{eq:scs2}
\end{align}
In other words, the realizability constraint, \Eqre{eq:realizable}, on
\Eqre{eq:sIto} mathematically corresponds to Eqs.\
(\ref{eq:scs1}--\ref{eq:scs2}). A diffusion process, governed by \Eqre{eq:sIto},
that satisfies Eqs.\ (\ref{eq:scs1}--\ref{eq:scs2}), ensures that the fractions
$Y$ and $1-Y$ satisfy $0 \le Y \le 1$, provided each event of the ensemble at
$t=0$ satisfies $0 \le Y \le 1$. The equal signs in the constraints on the drift
in Eqs.\ (\ref{eq:scs1}--\ref{eq:scs2}) allow for absorbing barriers at $Y = 0$
and $Y = 1$, respectively. The constraints on the diffusion term imply that
$B(Y,t)$ must either be nonlinear in $Y$ or $B(Y,t) \equiv 0$ for all $Y$. In
other words, since the diffusion term must be non-negative, required by
\Eqre{eq:sIto}, it can only be nonzero inside the allowed sample space if it is
also nonlinear.

\subsubsection{Realizable ternary processes: $\bv{N=3}$}
For $N\!=\!3$ variables, the unit-sum-constrained sample space and its boundary
are sketched in \Fige{fig:simplex}. In this case individual samples of the joint
probability, $F(Y_1,Y_2)$, are governed by the system,
\begin{align}
\mathrm{d}Y_1(t) & = A_1(Y_1,Y_2,t) \mathrm{d}t + b_{11}(Y_1,Y_2,t)
\mathrm{d}W_1(t) + b_{12}(Y_1,Y_2,t) \mathrm{d}W_2(t) \label{eq:2Ito1} \\
\mathrm{d}Y_2(t) & = A_2(Y_1,Y_2,t) \mathrm{d}t + b_{21}(Y_1,Y_2,t)
\mathrm{d}W_1(t) + b_{22}(Y_1,Y_2,t) \mathrm{d}W_2(t). \label{eq:2Ito2}
\end{align}
The allowed samples space is two-dimensional (a triangle) whose boundary,
defined by \Eqre{eq:boundary}, consists of the loop of lines,
\begin{equation}
(Y_1 = 0, \enskip Y_2 = 0, \enskip Y_1 + Y_2 = 1).
\end{equation}
For $N=3$, the state vector, governed by Eqs.\ (\ref{eq:2Ito1}--\ref{eq:2Ito2})
augmented by $Y_3=1-Y_1-Y_2$, stays inside the allowed region if
\begin{align}
A_1(Y_1\!=\!0,Y_2,t) \ge 0 \qquad & \mathrm{and} \qquad B_{11}(Y_1\!=\!0,Y_2,t)
= B_{12}(Y_1\!=\!0,Y_2,t) = 0, \label{eq:2scm1} \\
A_2(Y_1,Y_2\!=\!0,t) \ge 0 \qquad & \mathrm{and} \qquad B_{21}(Y_1,Y_2\!=\!0,t)
= B_{22}(Y_1,Y_2\!=\!0,t) = 0, \label{eq:2scm2} \\
A_1(Y_1\!+\!Y_2\!=\!1,t) \le 0 \qquad & \mathrm{and} \qquad
B_{11}(Y_1\!+\!Y_2\!=\!1,t) = B_{12}(Y_1\!+\!Y_2\!=\!1,t) = 0,
\label{eq:2scm3}\\
A_2(Y_1\!+\!Y_2\!=\!1,t) \le 0 \qquad & \mathrm{and} \qquad
B_{21}(Y_1\!+\!Y_2\!=\!1,t) = B_{22}(Y_1\!+\!Y_2\!=\!1,t) = 0. \label{eq:2scm4}
\end{align}
The realizability constraint, \Eqre{eq:realizable}, on the system of Eqs.\
(\ref{eq:2Ito1}--\ref{eq:2Ito2}) mathematically corresponds to Eqs.\
(\ref{eq:2scm1}--\ref{eq:2scm4}). The three fractions, $Y_1$, $Y_2$, and $Y_3 =
1 - Y_1 - Y_2$, governed by Eqs.\ (\ref{eq:YN},\ref{eq:2Ito1},\ref{eq:2Ito2}),
remain fractions of unity if their drift and diffusion terms satisfy Eqs.\
(\ref{eq:2scm1}--\ref{eq:2scm4}).  Naturally, an initial ensemble that
satisfies, $0 \le Y_1$, $0 \le Y_2$, and $Y_1+Y_2 \le 1$, is required. The
constraints on the diffusion terms in Eqs.\ (\ref{eq:2scm1}--\ref{eq:2scm4})
show that both $B_1$ and $B_2$ must either be nonlinear in $Y_1$ and $Y_2$,
respectively, or $B_1(Y_1,t) \equiv 0$ and $B_2(Y_2,t) \equiv 0$, for all $Y_1$
and $Y_2$, respectively. Furthermore, if one were to construct a process with
$A_1=A_1(Y_1)$, $B_{11}=B_{11}(Y_1)$, and $B_{12}=B_{12}(Y_1)$, then either
$A_2$ or $B_{22}$ must be a function of both $Y_1$ and $Y_2$ if $1-Y_1-Y_2\ge0$
is to be maintained, required by $Y_1+Y_2+Y_3=1$ with $Y_1\ge0$, $Y_2\ge0$,
$Y_3\ge0$. In other words, the unit-sum constraint couples at least 2 of the 3
fractions, governed by the system of Eqs.\
(\ref{eq:YN},\ref{eq:2Ito1},\ref{eq:2Ito2}).

\subsubsection{Realizable multi-variable processes: $\bv{N>2}$}
The multivariate It\^o diffusion process, equivalent to the Fokker-Planck
equation \Eqr{eq:FP}, is \cite{Gardiner_09}
\begin{equation}
\mathrm{d}Y_\alpha(t) = A_\alpha(\bv{Y},t)\mathrm{d}t + \sum_{\beta=1}^{N-1}
b_{\alpha\beta}(\bv{Y},t)\mathrm{d}W_\beta(t), \qquad \alpha =
1,\dots,N-1, \label{eq:mIto}
\end{equation}
with $B_{\alpha\beta} = \sum_{\gamma=1}^{N-1} b_{\alpha\gamma}b_{\gamma\beta}$
and the vector-valued Wiener process, $\mathrm{d}W_\beta(t)$, with mean
$\irmean{\mathrm{d}W_\beta} = 0$ and covariance $\irmean{\mathrm{d}W_\alpha
\mathrm{d}W_\beta} = \delta_{\alpha\beta} \mathrm{d}t$. Here
$\delta_{\alpha\beta}$ is Kronecker's delta. The sample space of allowed
realizations is now bounded by the non-overlapping hyperplanes, defined by
\Eqre{eq:boundary}.  The conditions, analogous to Eqs.\
(\ref{eq:scs1}--\ref{eq:scs2}) and Eqs.\ (\ref{eq:2Ito1}--\ref{eq:2Ito2}) that
ensure realizability for multiple variables are
\begin{align}
A_\alpha\left(Y_\alpha\!=\!0,Y_{\beta\ne\alpha},t\right) & \ge 0 \qquad
\mathrm{and} \qquad
B_{\alpha\beta}\left(Y_\alpha\!=\!0,Y_{\beta\ne\alpha},t\right) = 0,
\label{eq:scm1}\\
A_\alpha\left(\sum_{\alpha=1}^{N-1}Y_\alpha\!=\!1,t\right) & \le 0 \qquad
\mathrm{and} \qquad
B_{\alpha\beta}\left(\sum_{\alpha=1}^{N-1}Y_\alpha\!=\!1,t\right) = 0,
\label{eq:scm2}\\
&\qquad\qquad\qquad\qquad\qquad\quad \alpha,\beta=1,\dots,N-1.\nonumber
\end{align}
The realizability constraint, \Eqre{eq:realizable}, on the system of
\Eqres{eq:mIto} mathematically corresponds to Eqs.\
(\ref{eq:scm1}--\ref{eq:scm2}).  A diffusion process, governed by
\Eqre{eq:mIto}, that satisfies Eqs.\ (\ref{eq:scm1}--\ref{eq:scm2}) ensures that
the fractions $Y_\alpha$ satisfy $0 \le Y_\alpha \le 1$, $\alpha = 1,\dots,N$,
with $Y_N=1-\sum_{\beta=1}^{N-1} Y_\beta$, provided each event of the initial
ensemble at $t=0$ satisfies $0 \le Y_\alpha \le 1$. As before, the equal signs
in the constraints on the drifts in Eqs.\ (\ref{eq:scm1}--\ref{eq:scm2}) allow
for absorbing barriers at the boundaries. The constraints on the diffusion term
imply that for any $\alpha$, $B_{\alpha\beta}(Y_\alpha,Y_{\beta\ne\alpha},t)$
must either be nonlinear in $Y_\alpha$ or
$B_{\alpha\beta}(Y_\alpha,Y_{\beta\ne\alpha},t) \equiv 0$ for all $Y_\alpha$.
In other words, since the diffusion term must be non-negative semi-definite,
required by \Eqre{eq:mIto}, it can only be nonzero inside the allowed sample
space if it is also nonlinear. Eqs.\ (\ref{eq:scm1}--\ref{eq:scm2}) also show,
that while it is conceivable that $A_\alpha \ne A_\alpha(Y_\beta)$ and
$B_{\alpha\beta} \ne B_{\alpha\beta}(Y_\beta)$ for a single $\alpha$ and all
$\beta\ne\alpha$, if $\sum Y_\alpha=1$ is to be satisfied, either
$A_\alpha=A_\alpha(Y_\beta)$ or $B_{\alpha\beta} = B_{\alpha\beta}(Y_\beta)$
must hold for all $\beta\ne\alpha$.  In other words, the unit-sum constraint
couples at least $N-1$ equations of the system of Eqs.\ (\ref{eq:YN}) and
(\ref{eq:mIto}) governing $Y_\alpha$, $\alpha=1,\dots,N$.

Constraints on the functional forms of the drift and diffusion terms of the
multivariate Fokker-Planck equation \Eqr{eq:FP}, as a temporal representation
of random fractions, $Y_\alpha\!=\!Y_\alpha(t)$, have been developed. Eqs.\
(\ref{eq:scm1}--\ref{eq:scm2}) are our central result which ensure that sample
space events, generated by \Eqre{eq:FP} or its equivalent system of diffusion
processes, \Eqre{eq:mIto}, satisfy the realizability constraint at all times,
provided the initial ensemble is realizable. Since \Eqres{eq:FP} and
\Eqr{eq:mIto} govern $N-1$ variables and $Y_N = 1 - \sum_{\beta=1}^{N-1}
Y_\beta$, the unit-sum requirement, \Eqre{eq:unitsum}, is satisfied at all
times. An implication of Eqs.\ (\ref{eq:scm1}--\ref{eq:scm2}), exemplified in
\S\ref{sec:models}, is that random fractions represented by diffusion processes
must be coupled and nonlinear.

\section{Realizable evolution of statistics}
\label{sec:momentevolution}
Some implications of Eqs.\ (\ref{eq:scm1}--\ref{eq:scm2}) for the first few
statistical moments of the joint probability, governed by \Eqre{eq:FP}, are now
derived. This is useful for statistical moment equation representation of
fractions if individual samples and joint probabilities are not required.

\subsection{Realizable evolution of the means:
$\bv{\irmean{Y_\alpha}}$}
\label{sec:momenteqs1}
Multiplying \Eqre{eq:FP} by $Y_\gamma$ and integrating over all sample space,
see e.g.\ \cite{Pope_85}, yield the system of equations governing the means of
the fractions,
\begin{equation}
\frac{\partial\irmean{Y_\alpha}}{\partial t} = \rmean{A_\alpha} =
\mathcal{M}_\alpha, \qquad \alpha=1,\dots,N-1, \label{eq:means}
\end{equation}
where $A_\alpha\!=\!A_\alpha(\bv{Y},t)$. The evolution of the means can be made
consistent with the realizability constraint, \Eqre{eq:realizable}, if the means
are bounded and sum to one at all times. \Eqre{eq:means} shows that to keep the
means bounded, required by \Eqre{eq:boundedmeans}, the rate of change of the
means, $\mathcal{M}_\alpha$, must be governed by functions that satisfy
\begin{equation}
\lim_{\irmean{Y_\alpha}\to0} \mathcal{M}_\alpha = \lim_{\irmean{Y_\alpha}\to0}
\irmean{Y_\alpha},_t \ge 0 \qquad \textrm{and} \qquad
\lim_{\irmean{Y_\alpha}\to1} \mathcal{M}_\alpha = \lim_{\irmean{Y_\alpha}\to1}
\irmean{Y_\alpha},_t\le 0, \label{eq:constmeans}
\end{equation}
as the boundary of the state space is approached. In \Eqres{eq:constmeans}
$(\,\cdot\,),_t=\partial/\partial t$. \Eqres{eq:constmeans} imply that inside
the state space (i.e., away from the boundaries) $\mathcal{M}_\alpha$ must
either be a function of $\irmean{Y_\alpha}$ or $\mathcal{M}_\alpha \equiv 0$ for
all $t$. The means may also sum to one, required by \Eqre{eq:boundedmeans}, if
at least $N-2$ of \Eqres{eq:means} are coupled to each other. Consequently,
$\mathcal{M}_\alpha$ must be a function of $\irmean{Y_\beta}$ for all
$\beta\ne\alpha$.  \Eqre{eq:means} shows how the means are governed if a
Fokker-Planck equation \Eqr{eq:FP} or a diffusion process \Eqr{eq:mIto} governs
the underlying joint probability, e.g., only the mean of the drift, $A_\alpha$,
affects the evolution of the means.

\subsection{Realizable evolution of the second central moments:
$\bv{\irmean{y_\alpha y_\beta}}$}
\label{sec:momenteqs2}
Multiplying the Fokker-Planck equation \Eqr{eq:FP} by $y_\gamma y_\delta =
(Y_\gamma - \irmean{Y_\gamma})(Y_\delta - \irmean{Y_\delta})$ then integrating
over all sample space yields the equations governing the covariance matrix of
the fractions,
\begin{equation}
\frac{\partial\irmean{y_\alpha y_\beta}}{\partial t} = \rmean{y_\alpha A_\beta}
+ \rmean{y_\beta A_\alpha} + \rmean{B_{\alpha\beta}} =
\mathcal{C}_{\alpha\beta}, \qquad \alpha,\beta=1,\dots,N-1,
\label{eq:covariances}
\end{equation}
with $A_\alpha\!=\!A_\alpha(\bv{Y},t)$ and
$B_{\alpha\beta}\!=\!B_{\alpha\beta}(\bv{Y},t)$. The right hand side of
\Eqre{eq:covariances} is denoted by $\mathcal{C}_{\alpha\beta}$, the evolution
rate of the covariance matrix. \Eqre{eq:covariances} shows how the covariances
are governed if a Fokker-Planck equation \Eqr{eq:FP} or a diffusion process
\Eqr{eq:mIto} governs the underlying joint probability, e.g., $\irmean{y_\alpha
y_\beta}$ is symmetric at all times.  Following the development in
\S\ref{sec:realmoments}, four conditions must be satisfied by the system of
second moment equations \Eqr{eq:covariances} to ensure an evolution of the
covariances that is consistent with the realizability constraint,
\Eqre{eq:realizable}:
\begin{enumerate}
\item \emph{Symmetric covariance evolution.} The symmetry of the covariance
matrix can be ensured if $\irmean{y_\alpha y_\beta}(t=0)$ is symmetric, as well
as its evolution rates:
\begin{equation}
\mathcal{C}_{\alpha\beta} = \mathcal{C}_{\beta\alpha}.
\label{eq:constsymcovrates}
\end{equation}
\item \emph{Boundedness of the variances,} \Eqre{eq:boundedvariances}. This
condition can be ensured with
\begin{equation}
\lim_{\irv{y_\alpha}\to0} \mathcal{C}_{\alpha\alpha} = \lim_{\irv{y_\alpha}\to0}
\irv{y_\alpha},_t \ge 0 \qquad \textrm{and} \qquad
\lim_{\irv{y_\alpha}\to1} \mathcal{C}_{\alpha\alpha} = \lim_{\irv{y_\alpha}\to1}
\irv{y_\alpha},_t \le 0,
\label{eq:constvariances}
\end{equation}
as the boundary of the state space is approached, indicating that in general the
equation governing $\irv{y_\alpha}$ must either be a function of
$\irv{y_\alpha}$ or $\mathcal{C}_{\alpha\alpha}\equiv0$ for all $t$.
\item \emph{Boundedness of the covariances,} \Eqre{eq:boundedcovariances}. This
condition can be ensured if for $\alpha \ne \beta$,
\begin{equation}
\lim_{\irmean{y_\alpha y_\beta}\to-1} \mathcal{C}_{\alpha\beta} =
\lim_{\irmean{y_\alpha y_\beta}\to-1} \irmean{y_\alpha y_\beta},_t \ge 0 \qquad
\textrm{and} \qquad \lim_{\irmean{y_\alpha y_\beta}\to1}
\mathcal{C}_{\alpha\beta} = \lim_{\irmean{y_\alpha y_\beta}\to1}
\irmean{y_\alpha y_\beta},_t\le 0, \label{eq:constcovariances}
\end{equation}
as the boundary of the state space is approached, indicating that in general the
equation governing $\irmean{y_\alpha y_\beta}$ must either be a function of
$\irmean{y_\alpha y_\beta}$ or $\mathcal{C}_{\alpha\beta} \equiv 0$ for all $t$.
\item \emph{Zero row-sums,} \Eqre{eq:zerosums}. Differentiating
\Eqres{eq:zerosums} in time and using \Eqre{eq:covariances} yield the system
\begin{align}
\mathcal{C}_{11} + \mathcal{C}_{21} + \dots + \irmean{y_Ny_1},_t & = 0
\nonumber\\
\mathcal{C}_{12} + \mathcal{C}_{22} + \dots + \irmean{y_Ny_2},_t & = 0
\nonumber\\
& \vdots \label{eq:constzeroderivativesums}\\
\mathcal{C}_{1(N-1)} + \mathcal{C}_{2(N-1)} + \dots + \irmean{y_Ny_{N-1}},_t & =
0 \nonumber\\
\irmean{y_1y_N},_t + \irmean{y_2y_N},_t + \dots + \irv{y_N},_t & = 0. \nonumber
\end{align}
Performing the same substitutions on \Eqres{eq:constzeroderivativesums} that
resulted in \Eqre{eq:weakzerosums} we obtain the weaker constraint
\begin{equation}
\sum_{\alpha=1}^{N-1} \sum_{\beta=1}^{N-1} \mathcal{C}_{\alpha\beta} -
\irv{y_N},_t = 0.\label{eq:constweakzeroderivativesums}
\end{equation}
\end{enumerate}
We see that the trivial specification, $\mathcal{C}_{\alpha\beta} \equiv 0$,
satisfies all the above conditions, but also fixes the covariance matrix at its
initial state for all $t \ge t_0$, which is of limited applicability.

\subsection{Bounded evolution of the third central moments,
$\bv{\irs{y_\alpha}}$}
\label{sec:momenteqs3}
Multiplying the Fokker-Planck equation \Eqr{eq:FP} by $y_\gamma^3 = (Y_\gamma -
\irmean{Y_\gamma})^3$ then integrating yields the system governing the third
central moments, $\irs{y_\alpha}$, as
\begin{equation}
\frac{\partial\irs{y_\alpha}}{\partial t} = 3 \irmean{y_\alpha^2 A_\alpha} + 3
\sum_{\beta=1}^{N-1} \irmean{y_\alpha B_{\beta\beta}} = \mathcal{S}_\alpha,
\qquad \alpha = 1,\dots,N-1, \label{eq:thirdmoments}
\end{equation}
with $A_\alpha\!=\!A_\alpha(\bv{Y},t)$ and
$B_{\alpha\beta}\!=\!B_{\alpha\beta}(\bv{Y},t)$. The right hand sides of
\Eqres{eq:thirdmoments} are the evolution rates of the third moments, denoted by
$\mathcal{S}_\alpha$. The boundedness of the third moments, required by
\Eqre{eq:boundedthirdmoments}, can be ensured if
\begin{equation}
\lim_{\rs{y_\alpha}\to-1} \mathcal{S}_\alpha = \lim_{\rs{y_\alpha}\to-1}
\irs{y_\alpha},_t \ge 0 \qquad \textrm{and} \qquad
\lim_{\rs{y_\alpha}\to1} \mathcal{S}_\alpha = \lim_{\rs{y_\alpha}\to1}
\irs{y_\alpha},_t\le 0, \label{eq:constthirdmoments}
\end{equation}
as the boundary of the state space is approached, indicating that in general the
equation governing $\irs{y_\alpha}$ must either be a function of
$\irs{y_\alpha}$ or $\mathcal{S}_\alpha \equiv 0$ for all $t$. The conditions in
\Eqre{eq:constthirdmoments} only ensure boundedness, consequently, they are
necessary but not sufficient conditions for realizability of the third moments
as required by \Eqre{eq:realizable}. Note that the requirement on bounded sample
space has no implications on the boundedness of the skewness:
\begin{equation}
-\infty < \frac{\irs{y_\alpha}}{\irv{y_\alpha}^{3/2}} < \infty,
\label{eq:skewness}
\end{equation}
since $\irv{y_\alpha}\!\ge\!0$, \Eqre{eq:boundedvariances}.

\subsection{Bounded evolution of the fourth central moments,
$\bv{\irk{y_\alpha}}$}
\label{sec:momenteqs4}
Multiplying the Fokker-Planck equation \Eqr{eq:FP} by $y_\gamma^4 = (Y_\gamma -
\irmean{Y_\gamma})^4$ then integrating yields the system governing the fourth
central moments, $\irk{y_\alpha}$, as
\begin{equation}
\begin{split}
\frac{\partial\irk{y_\alpha}}{\partial t} = 4 \irmean{y_\alpha^3 A_\alpha} + 6
\sum_{\beta=1}^{N-1} \irmean{y_\alpha^2 B_{\beta\beta}} = \mathcal{K}_\alpha,
\qquad \alpha = 1,\dots,N-1, \label{eq:fourthmoments}
\end{split}
\end{equation}
with $A_\alpha\!=\!A_\alpha(\bv{Y},t)$ and
$B_{\alpha\beta}\!=\!B_{\alpha\beta}(\bv{Y},t)$.  The right hand sides of
\Eqres{eq:fourthmoments} are the evolution rates of the fourth moments, denoted
by $\mathcal{K}_\alpha$. The boundedness of the fourth moments, required by
\Eqre{eq:boundedneven}, can be ensured if
\begin{equation}
\lim_{\rk{y_\alpha}\to0} \mathcal{K}_\alpha = \lim_{\rk{y_\alpha}\to0}
\rk{y_\alpha},_t \ge 0 \qquad \textrm{and} \qquad
\lim_{\rk{y_\alpha}\to1} \mathcal{K}_\alpha = \lim_{\rk{y_\alpha}\to1}
\rk{y_\alpha},_t\le 0,
\label{eq:constfourthmoments}
\end{equation}
as the boundary of the state space is approached, indicating that in general the
equation governing $\irk{y_\alpha}$ must either be a function of
$\irk{y_\alpha}$ or $\mathcal{K}_\alpha \equiv 0$ for all $t$.  The conditions
in \Eqre{eq:constfourthmoments} only ensure boundedness, consequently, they are
necessary but not sufficient conditions for realizability of the fourth moments
as required by \Eqre{eq:realizable}. Note that, similar to the skewness in
\Eqre{eq:skewness}, the requirement on bounded sample space has no implications
on the upper bound of the kurtosis:
\begin{equation}
0 \le \frac{\irk{y_\alpha}}{\irv{y_\alpha}^{2}} < \infty,
\label{eq:kurtosis}
\end{equation}
since $\irv{y_\alpha}\!\ge\!0$, \Eqre{eq:boundedvariances}.

\subsection{Summary on realizable statistics of fractions}
The unit-sum constraint, \Eqre{eq:unitsum}, applied to a set of non-negative
random variables, bounds and constrains their statistical moments, as shown in
\S\ref{sec:realmoments}, as well as their time-evolutions. We examined the
evolution of the moments, $\irmean{Y_\alpha}$, $\irmean{y_\alpha y_\beta}$,
$\irs{y_\alpha}$, $\irk{y_\alpha}$, and showed how they are governed if an
underlying diffusion process is known.

Realizability of the means, as defined by \Eqre{eq:realizable}, can be ensured
if \Eqres{eq:boundedmeans} and \Eqr{eq:constmeans} are satisfied. Realizability
of the covariances can be ensured if Eqs.\
(\ref{eq:boundedvariances}--\ref{eq:zerosums}) and
(\ref{eq:constsymcovrates}--\ref{eq:constzeroderivativesums}) are satisfied.
Boundedness of the third moments is ensured by Eqs.\
(\ref{eq:boundedthirdmoments}) and \Eqr{eq:constthirdmoments}, while boundedness
of the fourth moments is ensured by Eqs.\ (\ref{eq:boundedneven}) and
\Eqr{eq:constfourthmoments}. The procedure outlined above can be continued to
derive additional constraints for consistency of the third, fourth, mixed, and
higher moments with the unit-sum constraint. The constraints reflect the coupled
and nonlinear nature of random fractions, both as instantaneous variables and
their statistics.

\section{A survey of realizable diffusion processes}
\label{sec:models}
A survey of existing diffusion processes that satisfy the realizability
constraints on drift and diffusion on the state-space boundary, Eqs.\
(\ref{eq:scm1}--\ref{eq:scm2}), is now given.

\subsection{Realizable binary process, $\bv{N\!=\!2}$: Beta}
An example for $N=2$, satisfying the realizability constraints on the drift and
diffusion terms on the sample-space boundary in Eqs.\
(\ref{eq:scs1}--\ref{eq:scs2}), is given in \cite{Bakosi_beta}, specifying the
drift and diffusion as
\begin{equation}
A(Y)=\frac{b}{2}(S-Y) \quad \textrm{and} \quad B(Y)=\kappa Y(1-Y),
\end{equation}
yielding the stochastic differential equation,
\begin{equation}
\mathrm{d}Y(t) = \frac{b}{2}(S-Y)\mathrm{d}t + \sqrt{\kappa Y(1-Y)}
\mathrm{d}W(t), \label{eq:beta}
\end{equation}
with $b>0$, $\kappa>0$, and $0<S<1$ excluding, while with $0 \le S \le 1$
allowing for absorbing barriers. In \Eqre{eq:beta} the drift is linear and the
diffusion is quadratic in $Y$. The invariant distribution of \Eqre{eq:beta} is
beta, which belongs to the family of Pearson distributions, discussed in detail by Forman \& S{\o}rensen \cite{Forman_08}. Of the special cases of the Pearson diffusions,
discussed in \cite{Forman_08}, only Case 6, equivalent to \Eqre{eq:beta},
produces realizable events. A symmetric variant of \Eqre{eq:beta}
was constructed in \cite{Cai_96}, which does not allow a non-zero skewness in
the statistically stationary state, see \cite{Bakosi_beta}.

\subsection{Realizable multivariate process, $\bv{N\!>\!2}$: Wright-Fisher}
A system of stochastic differential equations that satisfies the realizability
conditions for $N>2$ variables in Eqs.\ (\ref{eq:scm1}--\ref{eq:scm2}) is the
multivariate Wright-Fisher process \cite{Steinrucken_2013}, which specifies the
drift and diffusion terms as
\begin{equation}
A_\alpha(\bv{Y}) = \frac{1}{2}(\omega_\alpha - \omega Y_\alpha) \quad
\textrm{and} \quad B_{\alpha\beta}(\bv{Y}) = Y_\alpha(\delta_{\alpha\beta} -
Y_\beta),
\end{equation}
yielding the stochastic process,
\begin{equation}
\mathrm{d}Y_\alpha(t) = \frac{1}{2} (\omega_\alpha-\omega Y_\alpha) \mathrm{d}t
+ \sum_{\beta=1}^{N-1} \sqrt{Y_\alpha(\delta_{\alpha\beta}-Y_\beta)}
\mathrm{d}W_{\alpha\beta}(t), \qquad \alpha = 1,\dots,N-1,
\label{eq:WF}
\end{equation}
where $\omega = \sum_{\beta=1}^N \omega_\beta$, and $\omega_\beta>0$ are
parameters. \Eqre{eq:WF} is a generalization of \Eqre{eq:beta} for $N>2$
variables. The invariant distribution of \Eqre{eq:WF} is Dirichlet
\cite{Wright_49, Bakosi_dir}.

\subsection{Realizable multivariate process, $\bv{N\!>\!2}$: Dirichlet}
Another process that satisfies Eqs.\ (\ref{eq:scm1}--\ref{eq:scm2}), developed
in \cite{Bakosi_dir}, specifies the drift and diffusion terms as
{\allowdisplaybreaks
\begin{align}
A_\alpha(\bv{Y}) & = \frac{b_\alpha}{2}\big[S_\alpha Y_N -
(1-S_\alpha)Y_\alpha\big],\label{eq:a}\\\noalign{\smallskip}
B_{\alpha\beta}(\bv{Y}) & = \left\{
\begin{array}{lr}
\kappa_\alpha Y_\alpha Y_N & \quad \mathrm{for} \quad \alpha =
\beta,\\\noalign{\smallskip}
0 & \quad \mathrm{for} \quad \alpha \ne \beta,
\end{array}
\right.\label{eq:B}
\end{align}}%
resulting in the system of stochastic differential equations,
\begin{equation}
\mathrm{d}Y_\alpha(t) = \frac{b_\alpha}{2} \big[S_\alpha Y_N - (1-S_\alpha)
Y_\alpha\big] \mathrm{d}t + \sqrt{\kappa_\alpha Y_\alpha Y_N}
\mathrm{d}W_\alpha(t), \qquad \alpha=1,\dots,N-1, \label{eq:dir}
\end{equation}
with parameter vectors $b_\alpha > 0$, $\kappa_\alpha > 0$, and $0 < S_\alpha <
1$, and $Y_N$ given by \Eqre{eq:YN}. \Eqre{eq:dir} is also a generalization of
\Eqre{eq:beta} for $N>2$ variables. The invariant distribution of \Eqre{eq:dir}
is also Dirichlet, provided the parameters of the drift and diffusion terms
satisfy
\begin{equation}
(1-S_1) b_1 / \kappa_1 = \dots = (1-S_{N-1}) b_{N-1} / \kappa_{N-1}.
\label{eq:dirconst}
\end{equation}
Note that while there is no coupling among the parameters, $\omega_\alpha$, of
the drift and diffusion terms in the Wright-Fisher \Eqre{eq:WF}, the parameters,
$b_\alpha,$ $S_\alpha$, and $\kappa_\alpha$, of \Eqre{eq:dir} must be
constrained by \Eqre{eq:dirconst} to keep its invariant distribution Dirichlet.

\subsection{Realizable multivariate process, $\bv{N\!>\!2}$: Lochner's
generalized Dirichlet}
A generalization of \Eqre{eq:dir} is developed in \cite{Bakosi_gdir}, where the
drift and diffusion terms are given by
{\allowdisplaybreaks
\begin{align}
A_\alpha(\bv{Y}) & = \frac{\mathcal{U}_\alpha}{2}\left\{ b_\alpha\Big[S_\alpha
\mathcal{Y}_K - (1-S_\alpha)Y_\alpha\Big] + Y_\alpha\mathcal{Y}_K
\sum_{\beta=\alpha}^{K-1} \frac{c_{\alpha\beta}}{\mathcal{Y}_\beta}\right\},
\label{eq:ga}\\
B_{\alpha\beta}(\bv{Y}) & = \left\{ \begin{array}{lr}
\kappa_\alpha Y_\alpha \mathcal{Y}_K
\mathcal{U}_\alpha & \quad \mathrm{for} \quad \alpha = \beta,\\
\noalign{\smallskip} 0 & \quad
\mathrm{for} \quad \alpha \ne \beta,
\end{array} \right.\label{eq:gB}
\end{align}}%
with $\mathcal{Y}_\alpha = 1-\sum_{\beta=1}^\alpha Y_\beta$ and
$\mathcal{U}_\alpha = \prod_{\beta=1}^{K-\alpha} \mathcal{Y}_{K-\beta}^{-1}$,
yielding the stochastic process,
\begin{equation}
\begin{split}
\mathrm{d}Y_\alpha(t) = \frac{\mathcal{U}_\alpha}{2}\left\{
b_\alpha\Big[S_\alpha \mathcal{Y}_K -
(1-S_\alpha)Y_\alpha\Big] + Y_\alpha\mathcal{Y}_K
\sum_{\beta=\alpha}^{K-1}\frac{c_{\alpha\beta}}{\mathcal{Y}_\beta}\right\}\mathrm{d}t +
\sqrt{\kappa_\alpha Y_\alpha \mathcal{Y}_K
\mathcal{U}_\alpha}\mathrm{d}W_\alpha(t), \\
\qquad \alpha=1,\dots,K=N-1.
\label{eq:gendir}
\end{split}
\end{equation}
The invariant distribution of \Eqre{eq:gendir} is Lochner's generalized
Dirichlet distribution \cite{Lochner_75}, if the coefficients, $b_\alpha>0$,
$\kappa_\alpha>0$, $0<S_\alpha<1$, and $c_{\alpha\beta}$, with
$c_{\alpha\beta}=0$ for $\alpha>\beta$, $\alpha,\beta=1,\dots,K-1$, satisfy the
conditions developed in \cite{Bakosi_gdir}. Similar to \Eqre{eq:dir}, the
parameters of the drift and diffusion terms, $b_\alpha,$ $S_\alpha$,
$\kappa_\alpha$, and $c_{\alpha\beta}$, of \Eqre{eq:gendir} must be constrained
to keep the invariant distribution generalized Dirichlet. Setting
\begin{equation}
c_{1i}/\kappa_i=\dots=c_{ii}/\kappa_i=1 \quad\textrm{for}\quad
i=1,\dots,K\!-\!1,
\end{equation}
in \Eqre{eq:gendir} reduces to the standard Dirichlet process, \Eqre{eq:dir}.

All of Eqs.\ (\ref{eq:WF}), (\ref{eq:dir}), and (\ref{eq:gendir}) have coupled
and nonlinear diffusions terms. As discussed earlier, this is required to
simultaneously satisfy the realizability conditions in Eqs.\
(\ref{eq:scm1}--\ref{eq:scm2}), required to represent $N>2$ random fractions by
diffusion processes.

\section{Summary}
\label{sec:summary}
We have demonstrated that the problem of $N$ fluctuating variables constrained
by the unit-sum requirement can be discussed in a reduced sample space of $N-1$
dimensions. This allows working with the unique, universal, and mathematically
well-defined realizable sample space which produces samples and statistics
consistent with the underlying conservation principle.

We have studied multivariate diffusion processes governing a set of fluctuating
variables required to satisfy two constraints: (1) non-negativity and (2) a
conservation principle that requires the variables to sum to one, defined as
realizability. Our findings can be summarized as follows:
\begin{itemize}
\item The diffusion coefficients in stochastic diffusion processes, governing
fractions, must be coupled and nonlinear.
\item If the set of constraints
\begin{align}
A_\alpha\left(Y_\alpha\!=\!0,Y_{\beta\ne\alpha},t\right) & \ge 0 \qquad
\mathrm{and} \qquad
B_{\alpha\beta}\left(Y_\alpha\!=\!0,Y_{\beta\ne\alpha},t\right) = 0, \\
A_\alpha\left(\textstyle{\sum_{\alpha=1}^{N-1}}Y_\alpha\!=\!1,t\right) & \le 0
\qquad \mathrm{and} \qquad
B_{\alpha\beta}\left(\textstyle{\sum_{\alpha=1}^{N-1}}Y_\alpha\!=\!1,t\right) =
0, \\
&\qquad\qquad\qquad\qquad\qquad\quad \alpha,\beta=1,\dots,N-1,\nonumber
\end{align}
is satisfied as the state-space boundary is approached, the stochastic system,
\begin{equation}
\left\{ \begin{array}{rcl}
\mathrm{d}Y_\alpha(t)\!\!&\!\!=\!\!&\!\!A_\alpha(\bv{Y},t)\mathrm{d}t +
\sum_{\beta=1}^{N-1} b_{\alpha\beta}(\bv{Y},t)\mathrm{d}W_\beta(t), \qquad
\alpha = 1,\dots,N-1, \\[0.15cm]
Y_N\!\!&\!\!=\!\!&\!\!1-\sum_{\alpha=1}^{N-1}Y_\alpha,
\end{array} \right.
\end{equation}
with $B_{\alpha\beta}\!=\!\sum_{\gamma=1}^{N-1}
b_{\alpha\gamma}b_{\gamma\beta}$, ensures that the components of the vector of
fractions, $\bv{Y}\!=\!(Y_1,\dots,Y_N)$, remain non-negative and sum to one at
all times.
\item Boundedness of the sample space requires boundedness of the moments.
\end{itemize}
The constraints provide a method that can be used to develop drift and diffusion
functions for stochastic diffusion processes for variables satisfying a
conservation law and that are inherently realizable.

\bibliographystyle{unsrt}
\bibliography{jbakosi}

\end{document}